\documentclass[11pt,reqno]{amsart}
\usepackage{amsfonts,amsmath,amssymb}
\newtheorem{prop}{Proposition}

\def\beq#1#2\eeq{%
        \begin{equation}%
        \label{#1}%
            #2%
        \end{equation}%
    }

\usepackage{color}
\usepackage{graphicx}
\usepackage{epstopdf}

\theoremstyle{plain}
\newtheorem{theorem}{Theorem}

\theoremstyle{remark}

\theoremstyle{definition}

\title[Tropical Markov and Cayley]{Tropical Markov dynamics and Cayley cubic}

\author{ K. Spalding}\address{IFLScience Ltd, 6 Kean St, London WC2B 4AS, UK}
\email{kathrynspalding7@gmail.com}

\author{A.P. Veselov}
\address{Department of Mathematical Sciences,
Loughborough University, Loughborough LE11 3TU, UK  and Moscow State University, Moscow 119899, Russia}
\email{A.P.Veselov@lboro.ac.uk}

\begin{document}

\maketitle

\centerline{\it Dedicated to Emma Previato on her 65th birthday}

\begin{abstract}
We study the tropical version of Markov dynamics on the Cayley cubic, introduced by V.E. Adler and one of the authors. We show that this action is semi-conjugated to the standard action of $SL_2(\mathbb Z)$ on a torus, and thus is ergodic with the Lyapunov exponent and entropy given by the logarithm of the spectral radius of the corresponding matrix.
\end{abstract}

\section{Introduction}

In 1880 A.A. Markov \cite{M80} 
discovered a remarkable relation between the theory of binary quadratic forms and the following Diophantine equation known as the {\it Markov equation}
\begin{equation}
\label{Markov}
x^2 + y^2+ z^2 = 3xyz.
\end{equation}
Markov showed that all positive integer solutions can be found from the obvious one $x=y=z=1$ by applying the symmetry
\begin{equation}
\label{Markovinv}
(x, y, z) \rightarrow (x, y, 3xy-z)
\end{equation}
(which is a corollary of the Vieta formula for the Markov equation considered as a quadratic with respect to $z$) and permutations. The corresponding {\it Markov numbers}
\[1, 2, 5, 13, 29, 34, 89, 169, 194, 233, 433, 610, 985...\]
play a very important role in the theory of Diophantine approximations determining the rank of the ``most irrational" numbers (see for detail \cite{CF}). Many other relations were discovered later, including the theory of Frobenius manifolds and the related Painlev\'e-VI equation \cite{Dubr}, Teichm\"uller spaces \cite{F} and various problems in algebraic geometry  \cite{HP, Rud}.

The growth of Markov numbers was investigated by Don Zagier \cite{DZ}, who used the parallel (going back to Cohn \cite{C})
between the Markov tree and the Euclidean algorithm described by the equation
\begin{equation}
\label{Euclid}
a+b=c
\end{equation}
with coprime $a,b$.
One can view this parallel as a ``tropicalization" (known also as Maslov's ``dequantization" \cite{L}): if we write 
$$
x=e^{\frac{a}{\hbar}}, y=e^{\frac{b}{\hbar}}, z=e^{\frac{c}{\hbar}}
$$
and let $\hbar \to 0$, then we have from the Markov equation that
$$a+b=c$$
assuming that $a,b$ are less than $c$. Similar ideas were used by Andy Hone \cite{Hone}, who studied the growth problem in relation with Halburd's Diophantine approach to integrability \cite{Halb}.

In our paper \cite{SV} we used this relation to study the growth of the Markov numbers as functions of the paths on the Markov tree, where the Markov numbers are ``naturally growing", see Fig. 1. One can view this representation as a version of the Conway topograph \cite{Conway} for Markov triples.

\begin{figure}[h]
\begin{center}
 \includegraphics[height=38mm]{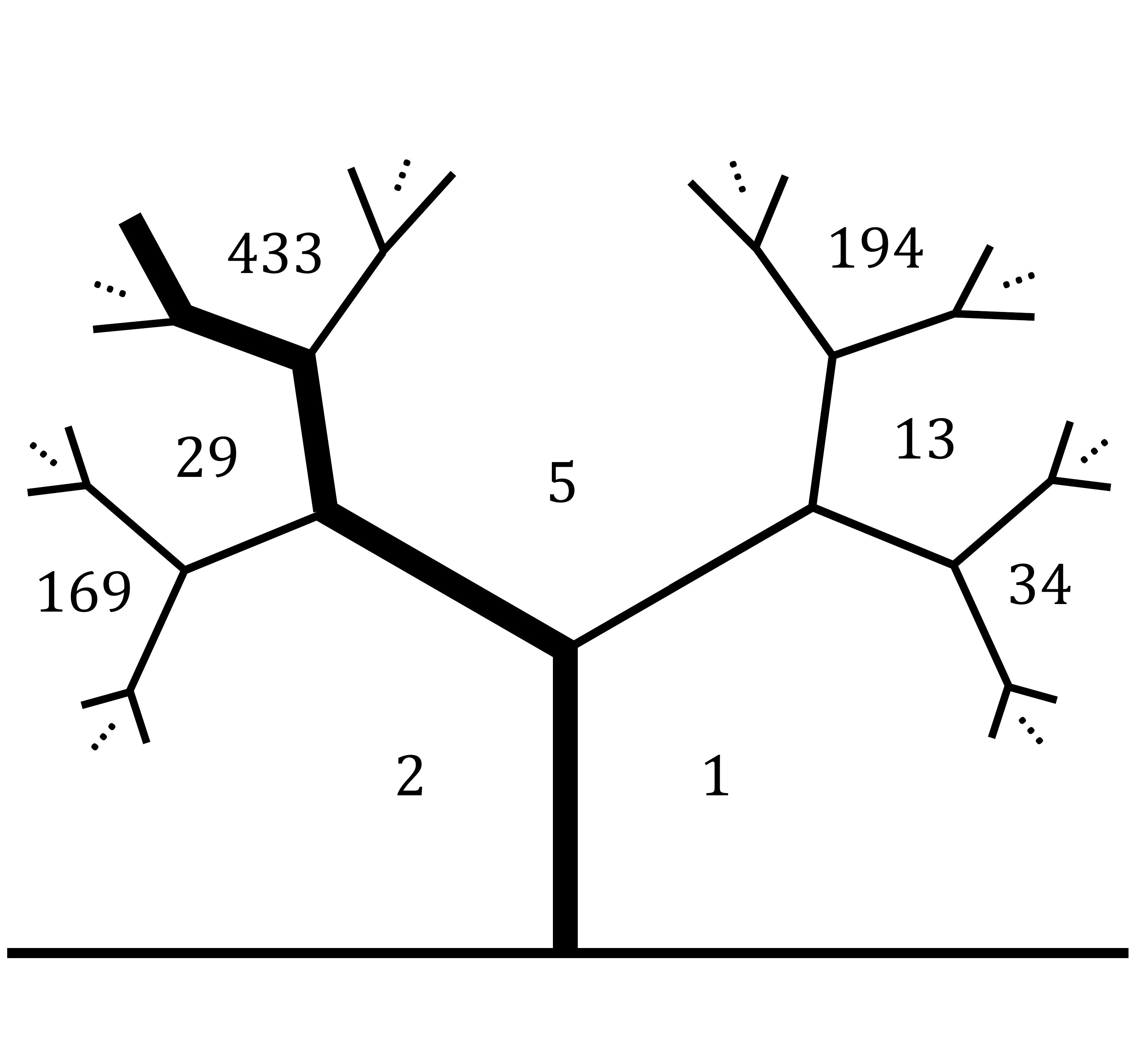}  \hspace{8pt}  \includegraphics[height=38mm]{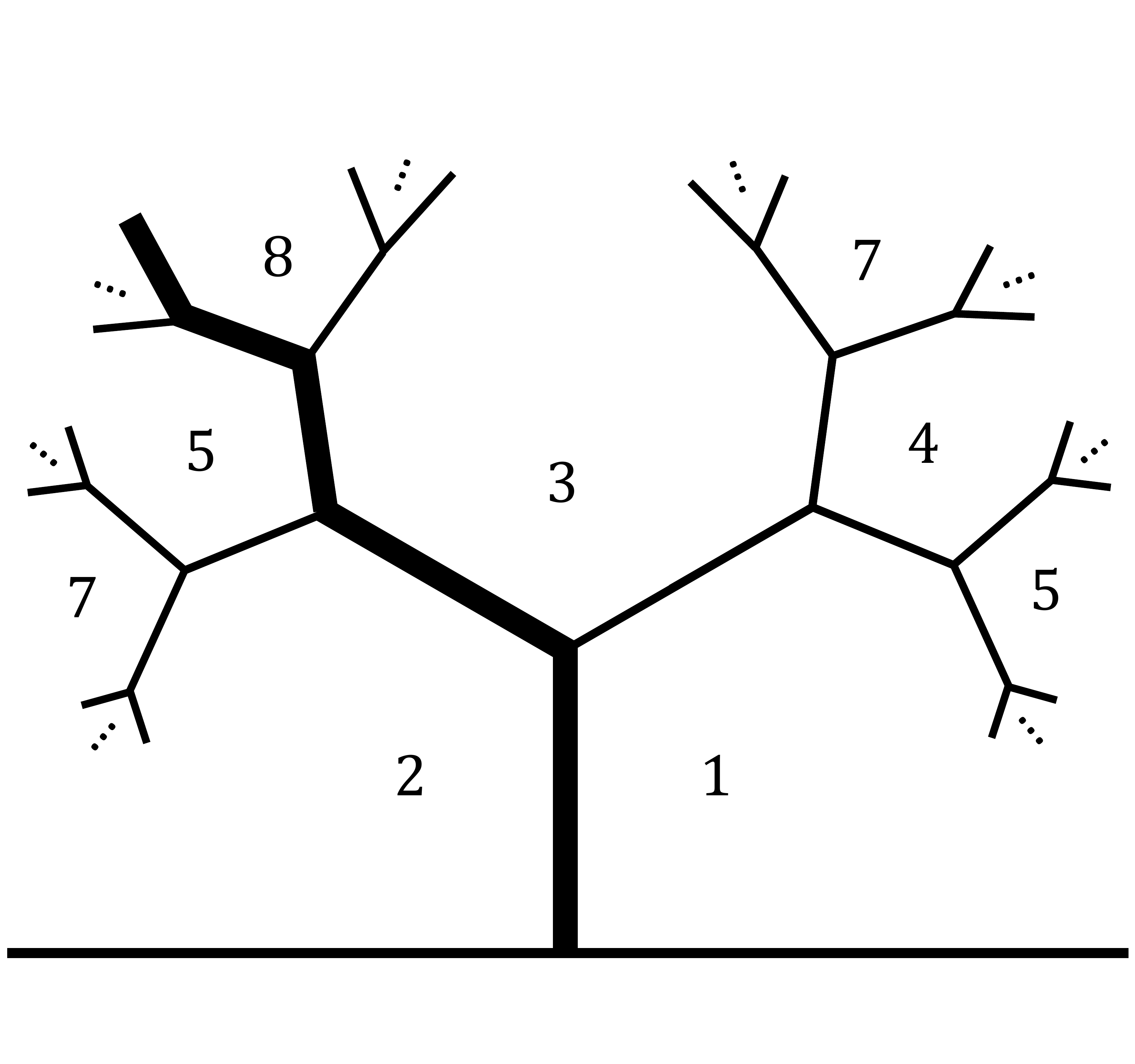}
\caption{\small Markov and Euclid trees with a path}
\end{center}
\end{figure}

More precisely, we defined the Lyapunov exponents of the Markov and Euclid trees
$\Lambda(\xi), \xi \in \mathbb RP^1$ as 
\beq{defL}
\Lambda(\xi)=\limsup_{n\to\infty}\frac{\ln(\ln z_n(\xi))}{n}=\limsup_{n\to\infty}\frac{\ln c_n(\xi)}{n}
\eeq
where $z_n(\xi), c_n(\xi)$ are the corresponding numbers along the path $\gamma_\xi$ on the Markov and Euclid trees respectively (see details in \cite{SV}).
The function $\Lambda(\xi)$ is $PGL_2(\mathbb Z)$-invariant and has some interesting properties studied in \cite{SV}.

In the present paper we consider the tropical version of the integrable case of Markov dynamics on the real surface given by
$$
x^2 + y^2+ z^2 = 3xyz+\frac{4}{9}
$$
or, equivalently after scaling the variables by 3:
\begin{equation}
\label{Cayley}
x^2 + y^2+ z^2 = xyz+4.
\end{equation}

From the algebro-geometric point of view the relation ({\ref{Cayley}) determines the classical surface known as the {\it Cayley cubic}. It was studied by Arthur Cayley in \cite{Cayley} and can be characterised as the cubic surface with 4 (which is maximal possible) conical singularities. The real version of Cayley cubic is shown in Fig. 2 prepared using MAPLE.

\begin{figure}[htbp]
\begin{center}
\includegraphics[width=8cm]{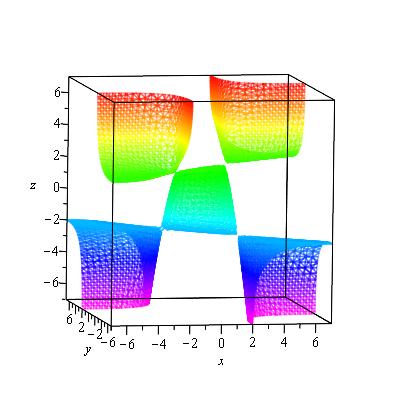}
\caption{Cayley cubic}
\end{center}
\end{figure}


It has four infinite sheets similar to the Markov surface (\ref{Markov}), but in this case the positive sheet (shown in the top right corner on Fig. 2) can be parametrized explicitly as
\begin{equation}
\label{coshp}
x=2\cosh a, \, y=2\cosh b, \, z=2\cosh c
\end{equation}
with $c=a+b$ (or in symmetric form $a+b+c=0$), which is in a good agreement with the tropical arguments above. This observation was used by Zagier \cite{DZ} to study the growth of Markov numbers and earlier by Mordell \cite{Mordell} for studying the Diophantine properties of this equation.

However, the Cayley cubic has also the middle part with the trigonometric parametrization
\begin{equation}
\label{Cayleypar}
x=2\cos a, \, y=2\cos b, \, z=2\cos c
\end{equation}
with $a+b+c=0.$ It bounds the body called the {\it spectrahedron}, which appears in a range of applications (see \cite{RG,Vinzant}).
It can be described as the set of semi-positive symmetric matrices of the form
$$
A=\begin{pmatrix}
  1 & x/2 & y/2 \\
  x/2 & 1 & z/2 \\
  y/2 & z/2 & 1 \\
\end{pmatrix}.
$$
Indeed one can check that the condition $\det A=0$ is equivalent to the Cayley relation (\ref{Cayley}). Note that $A$ with entries given by (\ref{Cayleypar}) is the Gram matrix for three unit vectors in $\mathbb R^3$ with pairwise angles $a,b,c$, so geometrically $\det A=0$ means that these vectors are coplanar and thus one of the angles is the sum of the others.

It is natural to ask if there is a ``tropical" analogue of the Markov dynamics on the middle part of the Cayley cubic. Here by the tropical dynamics we simply mean the dynamics determined by the piecewise linear maps.

Vsevolod Adler and one of the authors \cite{AV, V09} came up with a natural suggestion, simply replacing the spectrahedron by the regular tetrahedron with the vertices at the singular points, which are
$(2,2,2)$, $(2,-2,-2)$, $(-2,2,-2)$, $(-2,-2,2)$ (see Fig. 3). 

\begin{figure}[htbp]
\begin{center}
\includegraphics[width=6cm]{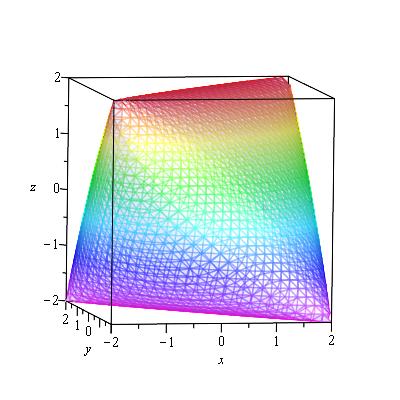}  \hspace{8pt}  \includegraphics[height=6cm]{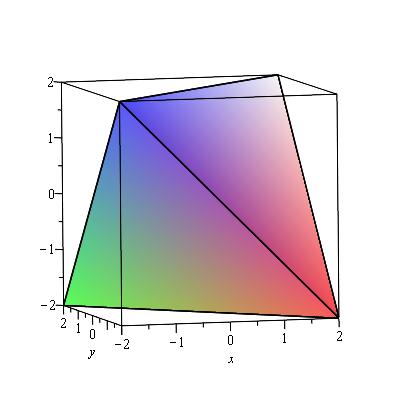}
\caption{Cayley's spectrahedron and tetrahedron}
\end{center}
\end{figure}

The corresponding boundary, which we will denote $T$, is determined by the ``tropical" Cayley equation
\begin{equation}
\label{tropCayley}
\max\{-u-v-w, -u+v+w, u-v+w, u+v-w\}=2.
\end{equation}
Note that this is {\it not} an analogue of the Cayley surface in the sense of tropical algebraic geometry \cite{IMS}, so the terminology might be a bit confusing (see more in Concluding remarks).

Following \cite{AV} consider the corresponding action of the modular group $PSL_2(\mathbb Z)$ generated by cyclic permutation of $u,v,w$ and the tropical Vieta involution
\begin{equation}
\label{tropinv}
(u,v,w) \rightarrow (u,v, -w+2f(u,v))
\end{equation}
where $f:\mathbb R^2 \rightarrow \mathbb R$ is a piecewise linear function defined by
\begin{equation}
\label{tropinvf}
 f(u,v)=\left\{\begin{array}{rrr}
   v & \text{if} &  u\ge|v|, \\
   u & \text{if} &  v\ge|u|, \\
  -v & \text{if} & -u\ge|v|, \\
  -u & \text{if} & -v\ge|u|.
 \end{array}\right.
\end{equation}

The plot of the function $f$ is shown in Fig. \ref{fig:f}.
\begin{figure}[h]
\centerline{\includegraphics[width=65mm]{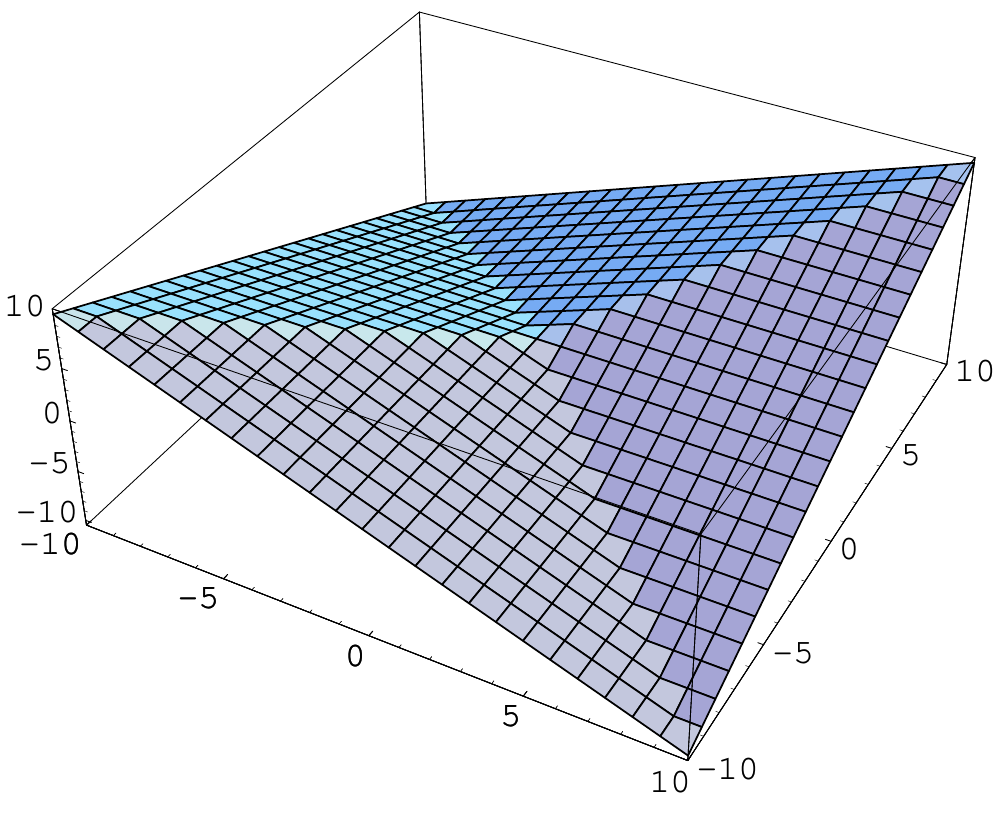}}
\caption{Function $f(x,y)$}
\label{fig:f}
\end{figure}

The aim of this paper is to study the properties of this action, which we will call {\it tropical Cayley-Markov dynamics.} 
Our main result is the following
\begin{theorem}
\label{main}
The tropical Cayley-Markov action of a hyperbolic element $A \in SL_2(\mathbb Z)$ on $T$ is ergodic, with the Lyapunov exponent and entropy given by the logarithm of the spectral radius of $A$. Their average growth along the path $\gamma_\xi$ on the planar binary tree is given by the function $\Lambda(\xi).$
\end{theorem}

The proof is by constructing the semi-conjugation of this action with the standard action of $SL_2(\mathbb Z)$ on a torus, using a natural tropical analogue of the parametrisation (\ref{Cayleypar}).

We should mention that the same idea was used by Cantat and Loray \cite{Can, CL} to compute the topological entropy of the (generalised) Markov dynamics (see also the important work of Iwasaki and Uehara \cite{I, IU} in this direction).

\section{Tropicalization of Markov dynamics and Cayley cubic}

Tropicalization (also known as dequantization \cite{L} or  ultra-discretization \cite{IKT}) can be applied to any dynamical system which can be written in algebraic form without a minus sign (subtraction-free), by replacing the operation of addition and multiplication by
$$X\oplus Y =\max (X,Y)$$
and 
$$X\otimes Y = X+Y$$
respectively. It is clear that this does not work directly for the Markov dynamics in the form (\ref{Markovinv}) because of the minus sign.

However one can consider another Vieta version (cf. Hone \cite{Hone}) 
\begin{equation}
\label{Markovinv2}
(x, y, z) \rightarrow (x, y, (x^2+y^2)/z),
\end{equation}
which can be naturally tropicalized as 
\begin{equation}
\label{Markovinvtrop}
(X, Y, Z) \rightarrow (X, Y, \max(2X, 2Y)-Z).
\end{equation}
Together with cyclic permutations of $X,Y,Z$, this generates the action of the modular group $PSL_2(\mathbb Z)$, which is known to be isomorphic to the free product $\mathbb Z_2*\mathbb Z_3.$ 


It has an invariant
$$
\Phi=\max(2X,2Y,2Z)-(X+Y+Z),
$$
or, equivalently,
\begin{equation}
\label{tropMint}\Phi= \max (X-Y-Z, Y-X-Z, Z-X-Y),
\end{equation}
which is the tropical version of the integral $$F=\frac{x^2+y^2+z^2}{xyz},$$
invariant under the Vieta involution (\ref{Markovinv}).

It is easy to see that the tropical equation $\Phi=0$ for positive integers $X,Y,Z$ defines the Euclidean algorithm and, as explained above, describes the asymptotic growth of the Markov triples in the logarithmic scale.

Let us now turn to the Cayley cubic case
$$
x^2+y^2+z^2=xyz+4.
$$
Adding 4 to the right hand side of the equation does not change much the asymptotic behaviour at infinity in the positive octant, and thus the tropicalization, which is the same as in the Markov case.
However, it changes the shape of the surface near the origin by adding the part bounded by 4 singular points (see Fig. 2 above).

Following \cite{AV}, replace this part by the surface $T$ of the tetrahedron with the same vertices.
The projection of $T$ to the $(u,v)$-coordinate plane is a 2-to-1 map to the corresponding square (see Fig. 5). 

\begin{figure}[h]
\centerline{\includegraphics[width=55mm]{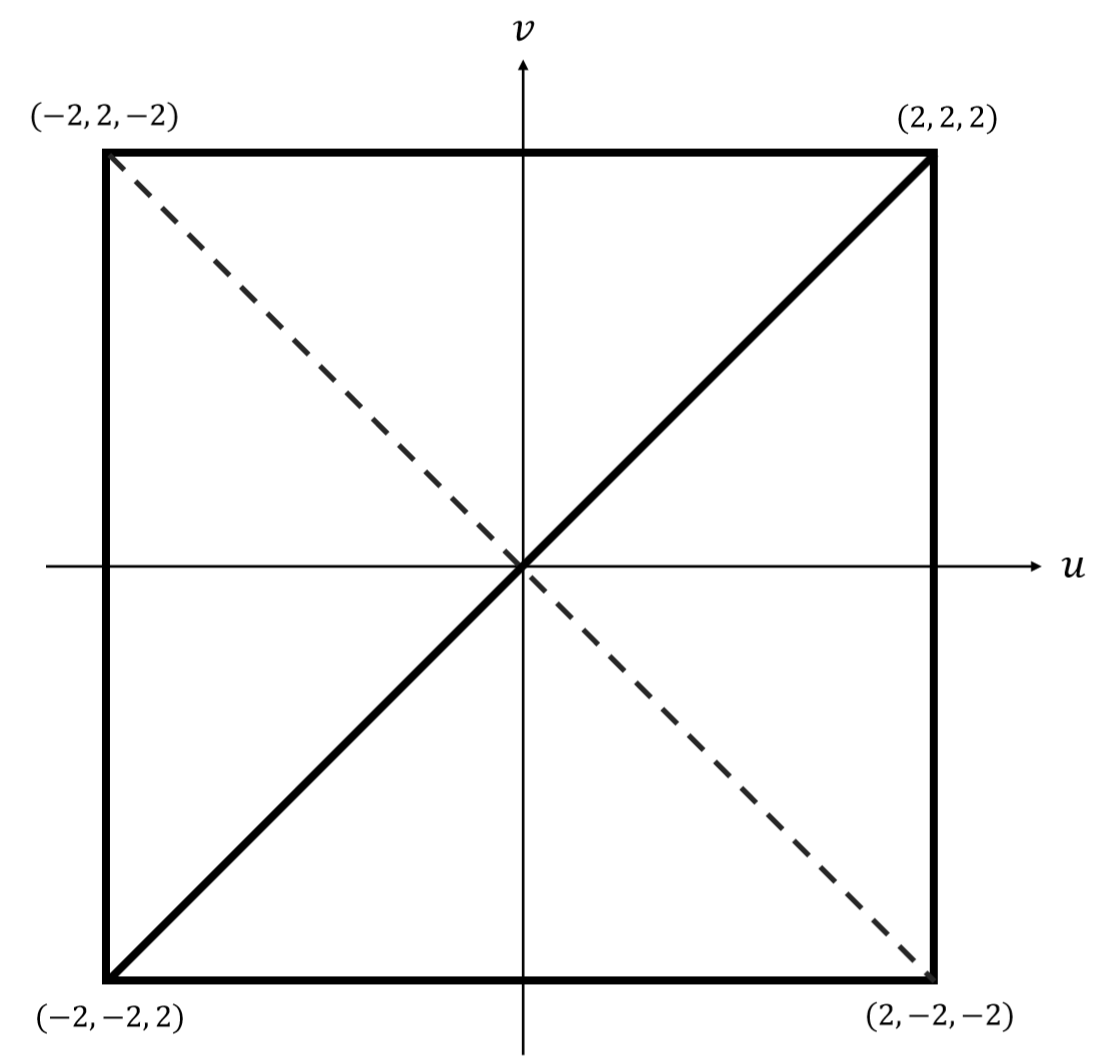}}
\caption{Projection of tropical Cayley surface $T$}
\end{figure}

One can check directly that the piecewise linear involution (\ref{tropinv}), (\ref{tropinvf})
swaps the branches of this double cover similarly to the Markov involution (\ref{Markovinv}), which was the motivation for introducing the tropical Cayley-Markov dynamics in \cite{AV}.

\begin{prop}
The function 
\begin{equation}
\label{tropCayleyint}
\Psi=\max\{ -u+v+w, u-v+w, u+v-w, -u-v-w\}
\end{equation}
is invariant under the tropical Cayley-Markov dynamics (\ref{tropinv}), (\ref{tropinvf}). The level set $\Psi=c, \, c>0$ is the surface of the regular tetrahedron with the vertices
$(c,c,c)$, $(c,-c,-c)$, $(-c,c,-c)$, $(-c,-c,c)$.
\end{prop}

The proof is by direct check. Note the difference of (\ref{tropCayleyint}) with (\ref{tropMint}), which can be rewritten equivalently as
$$
\Phi=-\min (-X+Y+Z, X-Y+Z, X+Y-Z).
$$


\section{Lyapunov exponents and entropy of the tropical Cayley-Markov dynamics.}

Now we would like to study the dynamical properties of the tropical Cayley-Markov action $PSL_2(\mathbb Z)=\mathbb Z_2*\mathbb Z_3$, where the action of $\mathbb Z_2$ is given by (\ref{tropinv}), (\ref{tropinvf}).

It is easy to see that this action preserves the usual Lebesgue measure on the surface of $T.$ 

The numerical calculations \cite{AV} showed the ergodic behaviour of the orbits of tropical Cayley-Markov dynamics at the level set $\Psi=c$, see Fig. 6.
\begin{figure}[t!]
\centerline{\includegraphics[width=80mm]{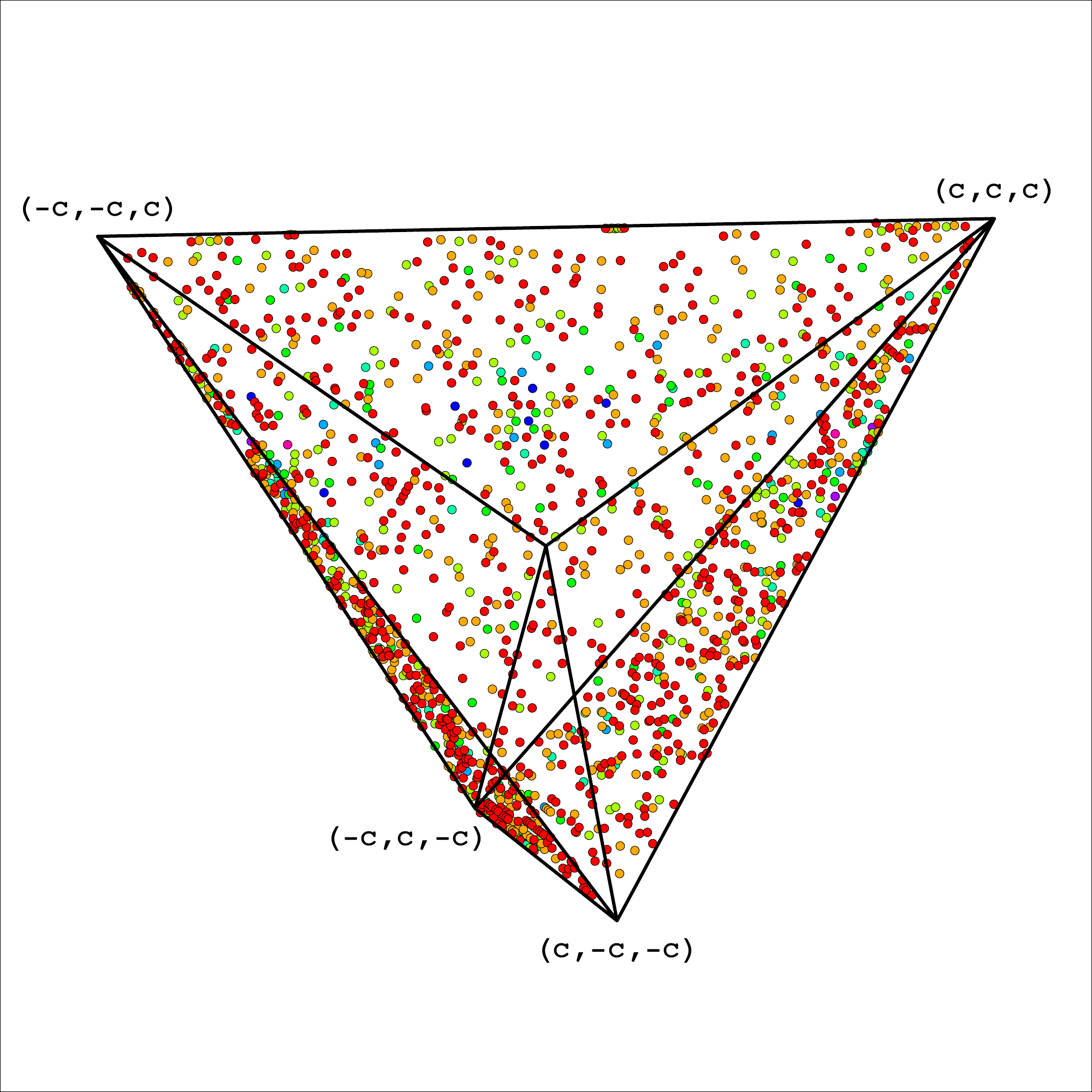}}
\caption{The level set of tropical Markov dynamics}
\label{fig:tetr}
\end{figure}

Now we are ready to prove this and our main Theorem \ref{main}. For this we need some results from our paper \cite{SV}.

Let us consider first the {\it Farey tree}, where at each vertex we have the fractions $\frac{p}{r}, \, \frac{q}{s}$ and their {\it Farey mediant}
$\frac{p+q}{r+s}$ (see Fig. 7).
Using the Farey tree we can identify the infinite paths $\gamma$ on a binary tree with real numbers $\xi \in [0,\infty]$ using the theory of continued fractions.
For example, for the golden ratio $\xi=\varphi:=\frac{\sqrt{5}+1}{2}$ we have the Fibonacci path shown in bold.

\begin{figure}[h]
\begin{center}
\includegraphics[height=48mm]{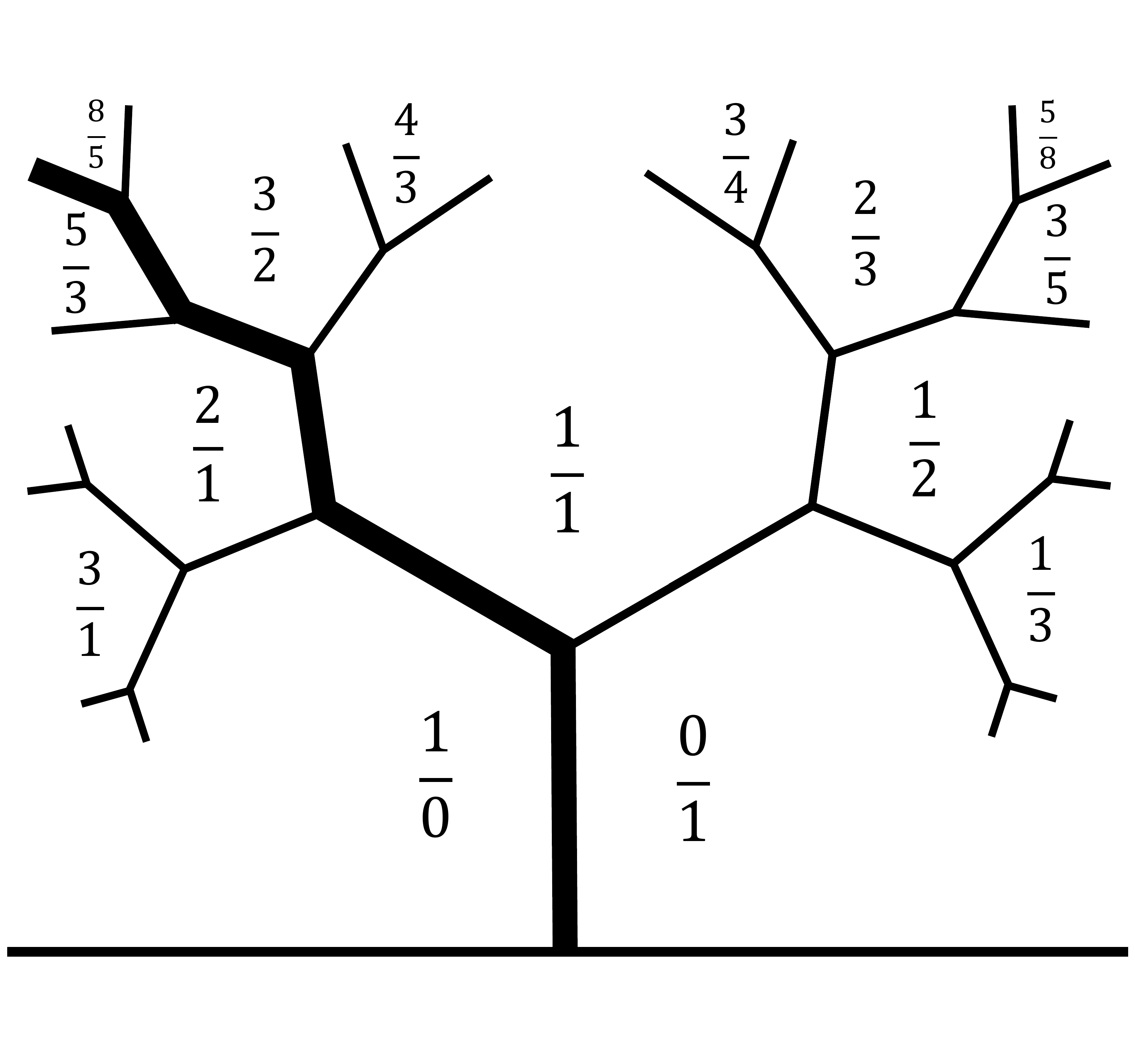}
\caption{\small The Farey tree with marked ``golden" Fibonacci path.}
\end{center}
\end{figure}

One can use the Farey tree to describe the monoid $SL_2(\mathbb N)$ consisting of matrices from $SL_2(\mathbb Z)$ with non-negative entries. 
Indeed, two neighbouring fractions $\frac{p}{r}, \, \frac{q}{s}$ can be combined into the matrix
\beq{A}
A= \begin{pmatrix}
  p & q \\
  r & s \\
\end{pmatrix} \in SL_2(\mathbb N).
\eeq

It can be shown \cite{SV} that the Lyapunov exponents (\ref{defL}) can be equivalently defined as
\beq{defeucl}
\Lambda(\xi)=\limsup_{n\to\infty}\frac{\ln \rho(A_n(\xi))}{n},
\eeq
where $A_n(\xi) \in SL_2(\mathbb N)$ is attached to the $n$-th edge along the path $\gamma_\xi$
and $\rho(A)$ is the {\it spectral radius} of the matrix $A$, defined as the maximum of the modulus of its eigenvalues.



Let's introduce the tropical cosine function $\cos_T  x$ as the period-2 piecewise linear even function given on the period by
$$
\cos_T x = 1-2|x|, \quad x \in [-1,1]
$$
 (see Fig. 8). It is known in Fourier theory as the {\it even triangle wave function}.

\begin{figure}[h]
\begin{center}
\includegraphics[height=30mm, width=100mm]{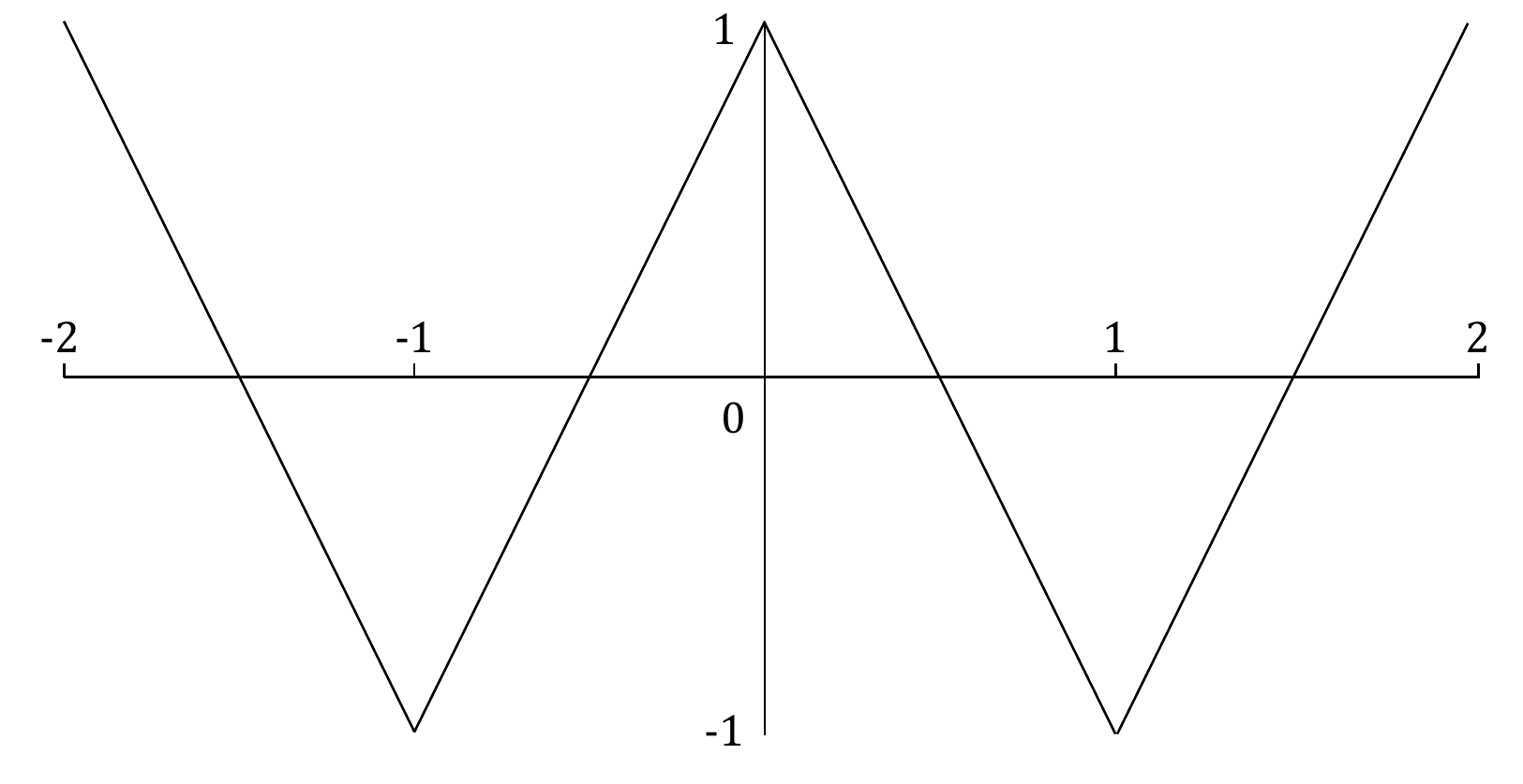} \hspace{20mm}
\caption{\small Tropical cosine (even triangle wave) function.}
\end{center}
\end{figure}
Define the tropical parametrization of $T$ by the following tropical analogue of (\ref{Cayleypar}):
\begin{equation}
\label{TCayleypar}
u=2 \cos_T \varphi, \, v=2 \cos_T\psi, \, w=2 \cos_T\chi, 
\end{equation}
where $\chi=\varphi+\psi$ and $(\phi, \psi) \in T^2=\mathbb R^2/(2\mathbb Z)^2.$

The corresponding map determines the 2-to-1 folding of the torus $T^2$ into the surface of the tetrahedron $T$ (see Fig. 9).

\begin{figure}[h]
\begin{center}
 \includegraphics[height=46mm]{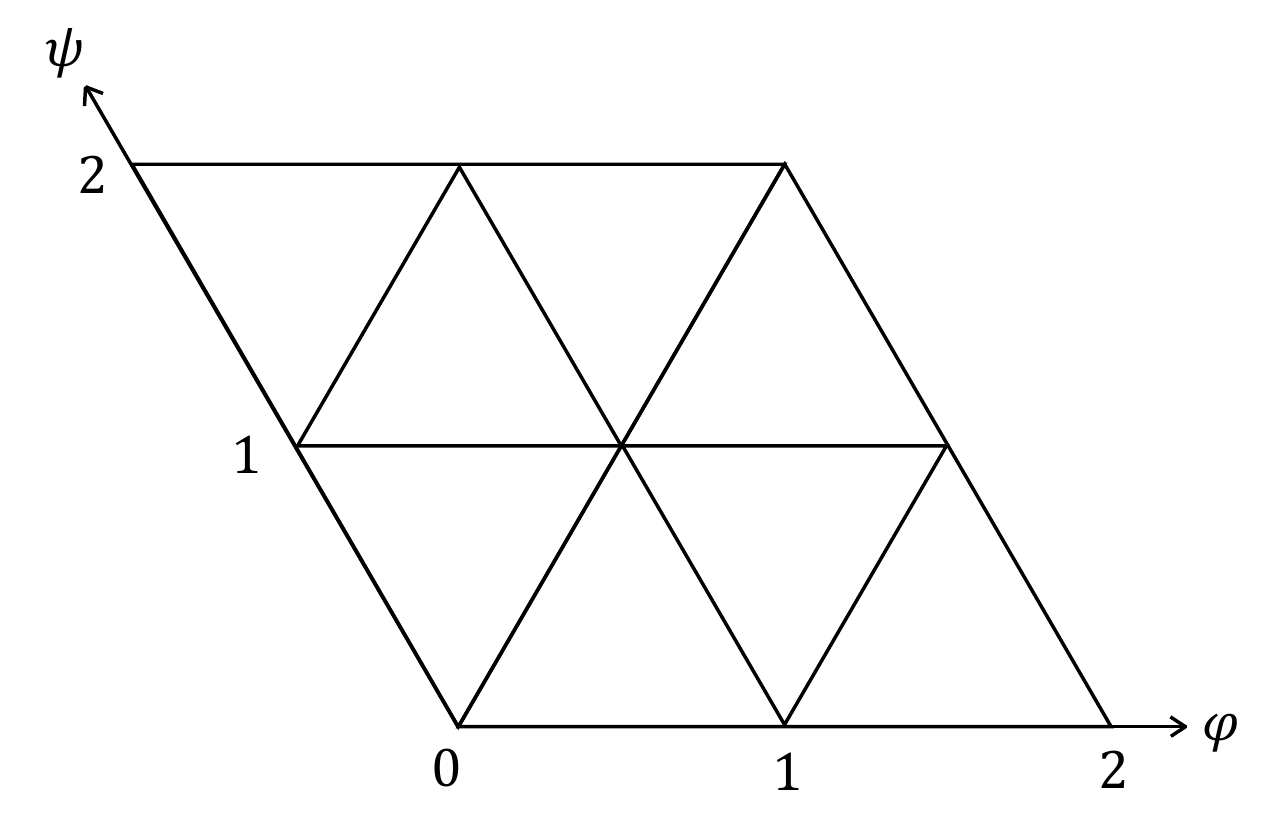}  \hspace{1pt}  \includegraphics[height=46mm]{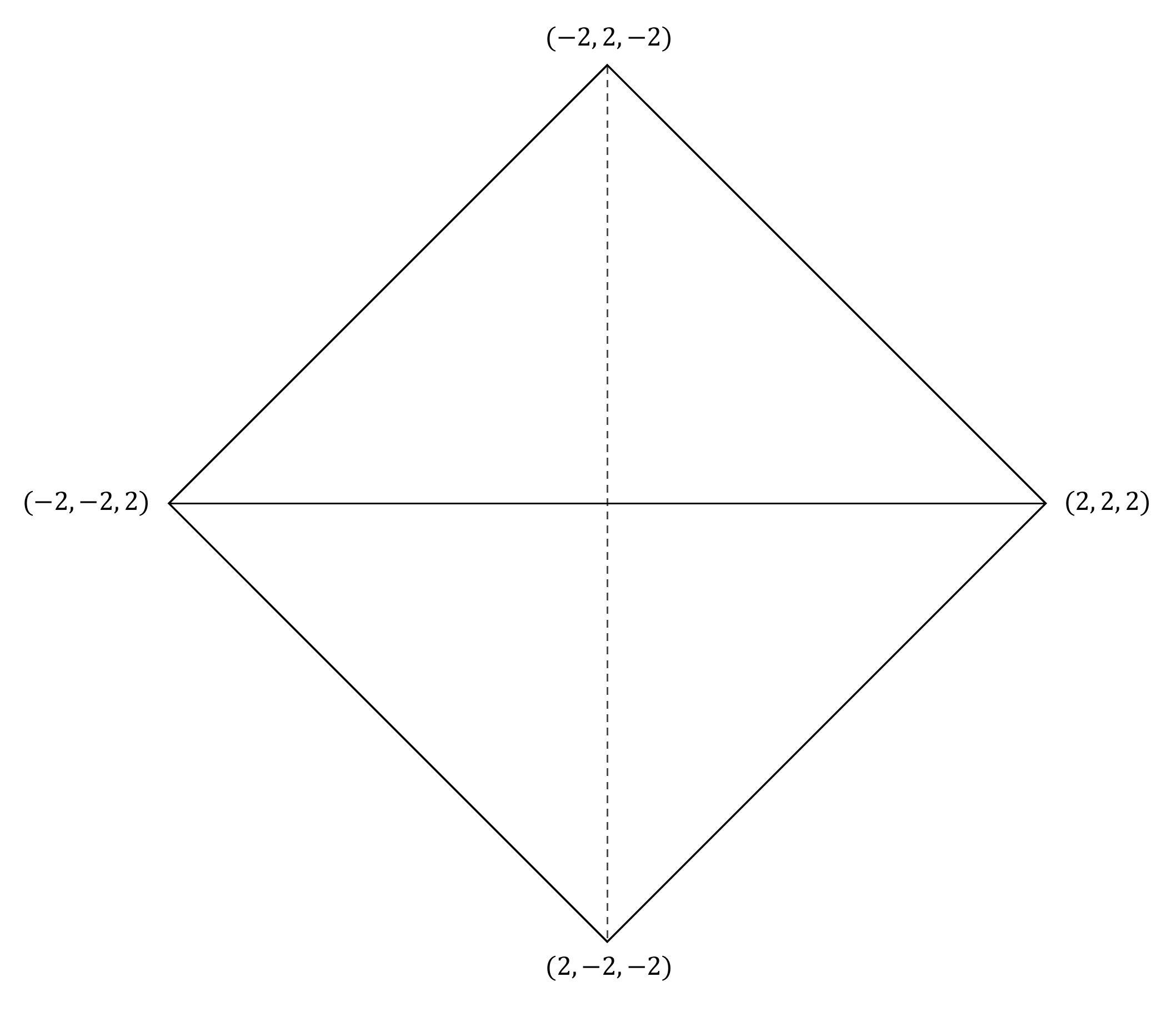}
\caption{\small Folding of the torus $T^2$ to the surface of tetrahedron $T$}
\end{center}
\end{figure}

The key observation now is the following

\begin{prop}
The parametrisation (\ref{TCayleypar}) semi-conjugates the tropical Cayley-Markov action of $A$ with the standard action of $A$ on the torus $T^2.$ 
\end{prop}

Indeed, $u$ and $v$ determine $\phi$ and $\psi$ by (\ref{TCayleypar}) uniquely up to a sign, which means that the two values of the corresponding coordinate $w$ are $w=2 \cos_T (\phi\pm\psi).$ Thus the tropical Cayley-Markov involution corresponds to the linear maps $(\pm \phi, \pm \psi, \pm(\phi+\psi)) \rightarrow (\pm \phi, \pm \psi, \pm(\phi-\psi))$, describing the action of $PSL_2(\mathbb Z)$ on the so-called lax superbases (see \cite{Conway}). 

Note that the surface of the tetrahedron $T$ is the quotient of the torus $T^2$ by the central symmetry group $\mathbb Z_2$, with fixed points corresponding to the vertices of the tetrahedron, so we have the following commutative diagram of the group actions
\begin{equation}
\label{commdia}
\begin{array}{ccc}
T^2&\stackrel{SL_2(\mathbb Z)}{\longrightarrow}&T^2\\ \downarrow
\lefteqn{\mathbb Z_2}& &\downarrow \lefteqn{\mathbb Z_2}\\
T&\stackrel{PSL_2(\mathbb Z)}{\longrightarrow}&T.
\end{array}
 \end{equation}

Since the action of a hyperbolic element $A \in SL_2(\mathbb Z)$ on a torus is known to be ergodic with the Lyapunov exponent and entropy given by the natural logarithm of the spectral radius of $A$ (see e.g. \cite{KH}), this completes the proof of Theorem 1.

\section{Concluding remarks}

A natural question is how exceptional is the case of Cayley-Markov dynamics considered in this paper.

Although replacing the spectrahedron by the regular tetrahedron looks like a natural tropicalization (in the sense of piecewise linearization), 
it cannot be explained by either traditional tropical algebraic geometry \cite{IMS} or the dequantisation/ultra-discretization procedure \cite{L,IKT}.

It would be nice therefore to see more similar examples in order to understand if there is any new general procedure here.

Although this is not directly related to Emma Previato's work, we believe that it is in the spirit of her very broad algebro-geometric view on integrable systems.
One of us (APV) had enjoyed many years of friendship with Emma, so we are very happy to dedicate this work for her 65th birthday and to wish her the best on this occasion.


\section{Acknowledgements} We are very grateful to Vsevolod Adler for many valuable discussions and crucial contribution to our preliminary work \cite{AV}, to Andy Hone and Masatoshi Noumi, who attracted our attention to the important papers \cite{Can, CL, I, IU}. 

This work was very much stimulated by the fruitful discussions with Leonid Chekhov, Boris Dubrovin, Andy Hone and Oleg Lisovyi. We are also grateful to Yiru Ye for help with preparation of the figures and to the referee for the helpful comments.

The work of K.S. was supported by the EPSRC as part of PhD study at Loughborough.

\end{document}